\documentclass[12pt]{article}
\usepackage{natbib,graphicx,setspace,amsmath,amssymb,subfigure,url,showlabels}

\newtheorem{thm}{Theorem}
\newtheorem{lem}{Lemma}
\newtheorem{rmk}{Remark}
\newtheorem{defi}{Definition}
\newtheorem{cor}{Corollary}
\newtheorem{prop}{Proposition}

\begin{document}

\title{Convergence of Nonparametric Functional Regression Estimates with Functional Responses }         
\author{
Heng Lian\\Division of Mathematical Sciences\\School of Physical and Mathematical Sciences\\Nanyang Technological University\\Singapore 637371\\Email: henglian@ntu.edu.sg\\Phone: +6565137175}        
\date{}          
\maketitle
\begin{abstract}
We consider nonparametric functional regression when both predictors and responses are functions. More specifically, we let $(X_1,Y_1),\ldots,(X_n,Y_n)$ be random elements in $\mathcal{F}\times\mathcal{H}$ where $\mathcal{F}$ is a semi-metric space and $\mathcal{H}$ is a separable Hilbert space. Based on a recently introduced notion of weak dependence for functional data, we showed the almost sure convergence rates of both the Nadaraya-Watson estimator and the nearest neighbor estimator, in a unified manner. Several factors, including functional nature of the responses, the assumptions on the functional variables using the Orlicz norm and the desired generality on weakly dependent data,  make the theoretical investigations more challenging and interesting.


\end{abstract}

\section{Introduction}
The problem of regression with functional predictors has been receiving increasing interests nowadays, boosted by more and more datasets with observations that can be naturally perceived as curves. This trend starts with the popular monograph \cite{ramsay02} that gives a detailed exposition of functional linear models. The existing literature contains numerous theoretical and empirical studies on functional linear models \citep{cardot99,cuevas02,james02,muller05,yao05,cai06,hall07,crambes09,bande10}. Nonparametric methods with functional predictors and scalar responses appear later \citep{ferraty02,ferraty04,ferraty06,preda07,biau10}, which by now has been widely accepted by the statistical community as a more flexible approach to functional regression with fewer structural assumptions imposed. As this area naturally develops and matures, the situation where the responses are also curves begins to receive more attention \citep{aguilera08,crambesmas09,hovath09}. For example, one might predict annual precipitation using temperature measurements as in \cite{ramsay05}, or predict future hourly electricity consumption based on past history as in \cite{antoch08}. Although these two studies follows the parametric approach to functional regression, it is clear that nonparametric approach is a viable alternative \citep{lian07}.

On the other hand, the assumption of independence in most theoretical investigations carried out so far is often too restrictive in many applications. The necessity to respond properly to data dependence is clearly demonstrated by the example given in \cite{ferraty02b} where a functional observation denotes the monthly electricity consumption over a year and thus it is unrealistic to assume that electricity consumption in one year is independent that of the previous year. In previous studies regarding nonparametric functional regression, dependence is incorporated based on some mixing conditions \citep{ferraty04}. Here we instead use the notion of $L^4-m-$approximability advocated in \cite{homann10,gabrys10} (with some appropriate minor extensions). The advantage compared to using mixing conditions is that the $L^4-m-$approximability condition is easily verified in many examples as shown in \cite{homann10}. 

In the more classical setting, the observation pairs reside in the Euclidean spaces. 
In this paper, we carry out a theoretical investigation of nonparametric functional regression with functional responses on dependent data. Two related classes of nonparametric estimates have been proposed, the k-nearest neighbor estimate (k-NN) and the Nadaraya-Watson kernel estimate. Because of their similarity in many aspects, we will try to unify the proofs for these two as much as possible. We will show almost sure convergence of these nonparametric estimators based on assumptions on Orlicz norms of the functional variables. Due to the functional nature of the responses and the assumption of weak dependence, the theoretical investigation poses serious challenges and some novel construction of martingale difference sequence will be introduced. Finally, we note that throughout the paper we use $C$ to denote a generic constant that assumes different values at different places.

\section{Almost sure convergence of nonparametric estimates}
\subsection{On the notion of Orlicz norm and weak dependence}
In this subsection we review the concept of Orlicz norm and collect some of its simple properties as a lemma here for easy reference later. Although all of the properties are simple and most are well-known, some others seem to be new (such as Lemma \ref{lem:psi} (vi)(vii)) which we cannot find in the existing literature. We also review and extend the notion of $L^4-m-$approximability of a data sequence using the more general Orlicz norm instead of $L^p$ norm.

Following \cite{vaartwellner96}, let $\psi$ be a convex, increasing function on $[0,\infty)$ with $\psi(0)=0$ and let $X$ be a real-valued random variable. The Orlicz norm (or $\psi$-Orlicz norm to emphasize its dependence on $\psi$) is defined as
\[\|X\|_\psi=\inf\{C>0:E[\psi(\|X\|/C)]\le 1\},\]
which can be shown to be indeed a norm. For random elements $X$ taking values in a normed space, the Orlicz norm of $\|X\|$ (which is a real-valued random variable) is also denoted by $\|X\|_\psi$ for simplicity. 

There are two commonly used $\psi$ function: $\psi(x)=x^p$ and $\psi(x)=\exp\{x^p\}-1$, $p\ge 1$, and throughout the paper we use $\psi_p$ to denote the latter. With $\psi(x)=x^p$, the Orlicz norm is simply the $L^p$ norm $(E[X^p])^{1/p}$. With $\psi(x)=\psi_p(x)=\exp\{x^p\}-1$, the finiteness of Orlicz norm of $X$ is closely related to the exponential decay of its tail probability, the exact statement of which is contained in the following Lemma together with other simple properties concerning Orlicz norm. 

\begin{lem}\label{lem:psi} Below we assume $\psi$ is a valid function that defines an Orlicz norm, that is, $\psi$ is convex, increasing on $[0,\infty)$ with $\psi(0)=0$. $X$ is a random variable.
\begin{enumerate}
\item[(i)] $P(|X|>x)\le 1/\psi(x/\|X\|_\psi), \forall x\ge 0$.
\item[(ii)] If $P(|X|>x)\le K \exp\{-Cx^p\}$ for all $x\ge 0$ and some constants $K$ and $C$, then $\|X\|_{\psi_p}\le ((1+K)/C)^{1/p}$.
\item[(iii)] If $\tilde\psi(x)=\psi(ax)$ for some $a>0$, then $\|X\|_{\tilde{\psi}}=a\|X\|_\psi$.
\item[(iv)] If $\tilde\psi(x)\le a\psi(x)$ for some $a\ge 1$, then $\|X\|_{\tilde{\psi}}\le a\|X\|_\psi$.
\item[(v)] If $\tilde{\psi}(x)=\phi(\psi(ax))$ for some $a>0$ and some concave increasing function $\phi$ with $\phi(0)=0$ and $\phi(1)=1$, then $\|X\|_{\tilde{\psi}}\le a\|X\|_\psi$.
\item[(vi)] If $\tilde{\psi}(x):=\psi(x^{1/p}), p\ge 1$ is convex, then $\left\||X|^p\right\|_{\tilde{\psi}}\le \|X\|_\psi^p$.
\item[(vii)] $\|E[X|\mathcal{G}]\|_\psi\le \|X\|_\psi$, for any $\sigma$-algebra $\mathcal{G}$.
\end{enumerate}
\end{lem}
\textit{Proof.} Results (i) and (ii) can be found in Section 2.2 of \cite{vaartwellner96}. (iii) is obvious by the definition of Orlicz norm. To prove (iv), we note that $E\tilde{\psi}(|X|/a\|X\|_\psi)\le aE\psi(|X|/a\|X\|_\psi)\le E\psi(|X|/\|X\|_\psi)\le 1$, where we used that $\psi(x/a)\le \psi(x)/a$ due to the convexity of $\psi$. For (v), 
since $E\tilde{\psi}(|X|/a\|X\|_\psi)=E\phi(\psi(|X|/\|X\|_\psi))\le\phi(E\psi(|X|/\|X\|_\psi))\le\phi(1)=1$ (using Jensen's inequality), we get $\|X\|_{\tilde{\psi}}\le a\|X\|_\psi$ by definition.
For (vi), the result follows from $E\tilde{\psi}(|X|^p/\|X\|_\psi^p)=E\psi(|X|/\|X\|_\psi)\le 1$. Finally, (vii) follows from $E\psi(E[X|\mathcal{G}]/\|X\|_\psi)=E\psi(E[X/\|X\|_\psi|\mathcal{G}])\le E(E(\psi(X/\|X\|_\psi)|\mathcal{G}))=E\psi(X/\|X\|_\psi)\le 1$, where we used $\psi(E[X/\|X\|_\psi|\mathcal{G}])\le E[\psi(X/\|X\|_\psi)|\mathcal{G}]$ due to convexity of $\psi$. $\Box$

We already noted that $L^p$ norm is a special case of Orlicz norm when $\psi(x)=x^p$. On the other hand, based on Lemma \ref{lem:psi} (v), one can show that $\|X\|_p\le C\|X\|_{\psi_q}$ for any $p,q\ge 1$ and $\|X\|_{\psi_{q_1}}\le C'\|X\|_{\psi_{q_2}}$ if $q_1\le q_2$, (where $C, C'$ are universal constants that only depends on $p,q,q_1,q_2$). In this sense the norm $\|.\|_{\psi_q}$ is stronger than $L^p$, and the more so with larger $q$.

As explained in the introduction, for data collected sequentially over time, the assumption of independence is not realistic. In \cite{homann10}, the authors formalize the notion of dependence for functional data using $L^4-m-$approximability. Instead of using the $L^4$ norm which is sufficient for the purpose of those studies, we instead use the Orlicz norm here.

\begin{defi}\label{def:app} Given a function $\psi$ that defines an Orlicz norm, a sequence $\{X_i\}_{i=1}^\infty$ (taking values in a normed space) with finite Orlicz norm is said to be $\psi-m-$approximable if we have the representation
\[X_i=h(\alpha_i,\alpha_{i-1},\ldots),\]
where the $\alpha_k$ are independent and identically distributed random elements of a measurable space and $h$ is a measurable function. In addition, we assume that if
\[X_i^{(m)}=h(\alpha_i,\alpha_{i-1},\ldots,\alpha_{i-m+1},\alpha_{i-m}',\alpha_{i-m-1}'\ldots),\]
with $\alpha_k'$ independent copies of $\alpha_0$, then
\[\sum_{m=1}^\infty \|X_m-X_m^{(m)}\|_\psi<\infty.\]

For a $\psi-m-$approximable sequence $\{X_i\}$, we say it is $\psi-m-$approximable with decay rate $\gamma_k$ if $\sum_{m=k}^\infty \|X_m-X_m^{(m)}\|_\psi=O(\gamma_k).$
\end{defi}

In \cite{homann10}, several examples of $L^p-m-$approximable sequence are given, minor modifications of these can produce more general $\psi-m-$approximable sequences. For example, a functional autoregressive process (Example 2.1 in \cite{homann10}) is $\psi-m-$approximable as long as the innovation noise has finite $\psi$-Orlicz norm, by the same arguments. Although not explicitly stated there, a functional autoregressive process is $\psi-m-$approximable with exponential decay rate: $\gamma_m=O(\exp\{-Cm\})$ for some constant $C$. 

\subsection{Nonparametric estimates and convergence rate}
Let $(X_1,Y_1),\ldots,(X_n,Y_n)$ be a stationary (in a strong sense) sequence of $\mathcal{F}\times\mathcal{H}$-valued random elements with $E\|Y\|<\infty$, where $\mathcal{F}$ is a semi-metric space with semi-metric $d(.,.)$ and $\mathcal{H}$ is a Hilbert space with norm $\|.\|$. The regression function is $r(x)=E(Y|X=x)$ and we can write $Y_i=r(X_i)+\epsilon_i$ where $\epsilon_i=Y_i-E(Y_i|X_i)\in \mathcal{H}$ are mean zero noises (in the sense of Bochner integral, see \cite{ledoux91}). \textit{In this subsection, we always consider probabilities and expectations conditional on $\{X_i\}$, in effect treating it as fixed}. The asymptotic results stated are thus conditional on predictors even though we do not state this explicitly in the following. The implications of random predictors are treated in the next subsection after we present the general convergence results in this subsection.

The regression function can be estimated by local weighting of responses 
\begin{equation}\label{eqn:rhat}
\hat{r}(x)=\sum_{i=1}^nW_{ni}(x)Y_i,
\end{equation}
where $(W_{n1}(x),\ldots,W_{nn}(x))$ is a probability vector of weights. Note that $W_{ni}(x)$ can be a function of all $X_k, k=1,\ldots, n$, instead of $X_i$ only, as is the case for k-NN estimates (see the examples below). Since in this paper we only investigate pointwise convergence at a fixed point $x$, we will use the notation $(W_{n1},\ldots, W_{nn})$ in the following for simplicity. 

We rank $(X_i,Y_i),i=1,\ldots,n$, based on increasing value of $d(X_i,x)$ (ties are broken by indices) and obtain a vector $(R_1,\ldots,R_n)$ such that $X_{R_i}$ is the $i$th nearest neighbor of $x$. Let $v_{ni}=W_{nR_i}$, we can write (\ref{eqn:rhat}) equivalently as 
\begin{equation}\label{eqn:rhat2}
\hat{r}(x)=\sum_{i=1}^n v_{ni}Y_{R_i}.
\end{equation}
Our consideration of weak dependence leads to extra complications in the proofs. If the observations are independent, then obviously $Y_{R_i}$ are also independent. However, if $(Y_1,Y_2,\ldots)$ is merely stationary, then $(Y_{R_1},Y_{R_2},\ldots)$ is no longer stationary in general since the order of observations are broken. We will thus use representation (\ref{eqn:rhat}) in most parts of our proofs, although representation (\ref{eqn:rhat2}) is easier to manipulate in the study of k-NN estimates for independent data.

\textit{Example 1. Simple nearest neighbor estimate.} Take $v_{ni}=1/k$ for $i\le k$ and $v_{ni}=0$ for $i>k$, so that the regression function estimate is just the average of responses corresponding to the $k$ nearest neighbors of $x$. Even in this simplest case, although $v_{ni}$ is only a deterministic sequence, $W_{ni}$ still depends on all $X_j, 1\le j\le n$ since all predictors jointly determine $x$'s neighbors. More generally, we can take $v_{ni}$ to be a deterministic sequence with $v_{n1}\ge v_{n2}\ge\cdots\ge v_{nn}$ thus putting more weights on data closer to $x$.

\textit{Example 2. Nearest neighbor estimate based on kernel.} Take \\
$W_{ni}=K(d(X_i,x)/H)/\sum_j K(d(X_j,x)/H)$ where $K$ is a kernel function and $H$ is the distance of the $k$th nearest neighbor to $x$. Mathematically,
\begin{equation}\label{eqn:H}
H=\inf\{h\in R: \sum_{i=1}^nI\{X_i\in B(x,h)\}\ge k\},
\end{equation}
where $B(x,h)=\{x'\in\mathcal{F}: d(x',x)\le h\}$ and $I\{.\}$ denotes the indicator function. In this subsection, since we condition on predictors $\{X_i\}$, $H$ is a known fixed value.

\textit{Example 3. Nadaraya-Watson estimate.} Take $W_{ni}=K(d(X_i,x)/H)/\sum_j K(d(X_j,x)/H)$, which has exactly the same form as in the previous example. However, here $H$ is a predetermined value usually called the bandwidth parameter, not derived from distance of $x$'s $k$th nearest neighbor. Typically, one applies the same value of $H$ for all values of $x$. Thus compared to nearest neighbor estimate, the Nadaraya-Watson estimate is not adaptive to the local sparseness of data. In this subsection when conditioning on predictors and for a given $x$, of course Nadaraya-Watson estimator is same as that in Example 2 since $H$ is fixed in both cases. The differences will appear in the next subsection.

Naturally we need the following assumption on the regression function to obtain nontrivial rates of convergence.

\textit{Assumption 1: $r$ is bounded and Lipschitz continuous. That is $\|r(x)\|\le B, \forall x\in\mathcal{F}$ and $\|r(x)-r(x')\|\le M\|x-x'\|^\alpha$. }

In fact, since we only consider pointwise convergence, it suffices that $r$ is Lipschitz continuous on an open neighborhood of $x$. We nevertheless use the above assumption for simplicity in statements.

\textit{Assumption 2: We assume $v_{n1}\ge v_{n2}\ge\cdots\ge v_{nn}$. Moreover, for some integer $k$ with $k/n\rightarrow 0$ and $k/\log n\rightarrow \infty$, we have $\sum_{i=k+1}^n v_{ni}=O(b_n)$ and $(\sum_{i=1}^nv_{ni}^2)^{1/2}=O(c_{n2})$ with $b_n, c_{n2}\rightarrow 0$. Also, we denote by $H$ the distance to $x$ from its $k$th nearest neighbor, and we assume $H\rightarrow 0$. }

Although Assumption 2 as stated is more amenable for use for k-NN estimates, it can also be used for Nadaraya-Watson estimate, which will be clear in the next subsection. We also impose the following assumptions on the noise.

\textit{Assumption 3: Given a convex increasing function $\psi$ with $\psi(0)=0$, and suppose for some constants $C>0$, some concave increasing function $\phi$ with $\phi(0)=0, \phi(1)=1$, we have that $x^r\le \phi(\psi(Cx))$ for some $r\ge 2$. Moreover, $M:=\|\epsilon_i\|_\psi<\infty$ and the stationary sequence $(\epsilon_1,\epsilon_2,\ldots)$ is $\psi-m-$approximable with decay rate $\{\gamma_k\}$.}

In the above assumption, the Orlicz norm is used for bounding the tail probability of noises (Lemma \ref{lem:psi} (i)) as well as controlling the dependence. It is possible of course to use different $\psi$ for these two different purposes, but using the same $\psi$ seems to be most natural since it concern the same noises. The assumption $x^r\le \phi(\psi(Cx))$ deserves some explanation. By Lemma \ref{lem:psi} (v), this implies that the $r$-th moment of the noise variable is finite, for some $r\ge 2$ and it is in particular satisfied by $\psi(x)=x^p$ for $p\ge r$ and $\psi(x)=\psi_q(x)$ for $q\ge 1$. When a stronger $\psi$-Orlicz norm is used, Assumption 3 imposes a stronger constraint, but the summability conditions in Theorem \ref{thm:main} below are easier to satisfy.

Our main result for functional nonparametric estimates with functional responses is the following.
\begin{thm}\label{thm:main}
If Assumption 1, 2 and 3 hold and if one can find sequences $a_n\rightarrow 0, L_n\rightarrow 0, x_n\rightarrow 0, m_n \mbox{ with } m_n \mbox{ an integer between } 1 \mbox{ and } n $, such that (in the rest of the paper these sequences are simply denoted by $a,L,x,m$)
\begin{enumerate}
\item[(*)] The four sequences, $\exp\{-Ca^2/(aL+m^2c_{n2}^2+x)\}$ for some constant $C$ big enough, $1/\psi(\sqrt{x/2}/(\gamma_1c_{n2}))$, $(m/a)/\psi^{1-1/r}\left(L/(2Mmv_{n1})\right)$, and $1/\psi(a/(2nv_{n1}\gamma_{m}))$, are all summable over $n$. 
\end{enumerate}
Then $\|\hat{r}(x)-r(x)\|=O(b_n+H^\alpha+a+(\gamma_1v_{n1})^{1/2})$ almost surely.
\end{thm}

\begin{rmk}\label{rmk:1} Here we present a unified result for both nearest-neighbor estimate and the Nadaraya-Watson estimate. For nearest-neighbor estimate, $k$ is a pre-specified constant and typically $b_n$ and $c_{n2}$ are explicit functions of $k$ and thus deterministic. On the other hand, $H$ depends on $k$ through (\ref{eqn:H}) and thus depends on predictors. The situation for the Nadaraya-Watson estimate is exactly the opposite. $H$ will be prespecified (typically as a function of sample size) and $k$ is the number of predictors falling into the ball with radius $H$ and thus depends on data.  Similarly, $v_{ni}$ as order statistics of $W_{ni}$ depend on predictor values.
\end{rmk}

\begin{rmk} Because of the requirement $\sum_{n=1}^\infty \exp\{-Ca^2/(aL+m^2c_{n2}^2+x)\}<\infty$, we see that the sequence $a$ cannot converge faster than $mc_{n2}$ and thus we will focus on cases where this rate is achievable up to some logarithmic terms in the following.
\end{rmk}

\begin{rmk} For independent data, $\gamma_1=0$ and the term $(\gamma_1v_{n1})^{1/2}$ does not appear. More generally, this term can be ignored as long as $v_{n1}=O(c_{n2}^2)$, by Remark 1 above. As an example, we obviously have $v_{n1}=c_{n2}^2=1/k$ for the simplest k-NN estimate with $v_{ni}=1/k, i\le k$. In the next subsection, one will see that for the Nadaraya-Watson estimate in Example 3 above, we also have that $v_{n1}$ and $c_{n2}^2$ are of the same order under mild assumptions.
\end{rmk}

\begin{rmk} In the convergence rate, $b_n$ and $H^\alpha$ represent the bias while $a$ comes from the variance of the estimator. As presented above, which aims for generality rather than clarity, it is hard to see what the convergence rate is in typical situations, and thus we discuss the rates in some special cases in the rest of this subsection. 
\end{rmk}

\textbf{Independent case} When the data are independent, $1/\psi(\sqrt{x_n/2}/(\gamma_1c_{n2}))$ and $1/\psi(a/(2nv_{n1}\gamma_{m}))$ are zero (Informally, $\gamma_m=0$ when data are independent and we take $\psi(\infty)=\infty$. More rigorously, it can be seen from the proofs that these two terms are zero), and we can take $m=1$, $x=0$. Taking $L=c_{n2}$ and $a=(\log n)c_{n2}$, the first sequence in (*) is then obviously summable. So as long as $1/\left(a\psi^{1-1/r}\left(c_{n2}/(2Mv_{n1})\right)\right)$ is summable, we have convergence rate $(\log n)c_{n2}$. For the simplest nearest neighbor estimate with $v_{ni}=1/k, i\le k$, we have $c_{n2}=1/\sqrt{k}$. The expression $1/a/\psi^{1-1/r}\left(c_{n2}/(2Mv_{n1})\right)$ is simplified to $\sqrt{k}/\left((\log n)\psi^{1-1/r}(\sqrt{k}/2M)\right)$. For $\psi(x)=x^p$ or $\psi(x)=\exp\{x^p\}-1$, this obviously is a restriction on $k$, in particular that $k$ should diverge fast enough at a certain rate. We note that by existing results on k-NN estimate for independent data with scalar responses, the variance term is expected to be $c_{n2}=1/\sqrt{k}$, 
 which agrees with the rate here up to a logarithmic term. In summary, we have

\begin{cor}\label{cor:knn} For simplest k-NN estimate with $v_{ni}=1/k, i\le k$, if $\sum_{n=1}^n\sqrt{k}/\psi^{1-1/r}(\sqrt{k}/2M)<\infty$ where $M=\|\epsilon_i\|_\psi$, then $\|\hat{r}(x)-r(x)\|=O(H^\alpha+(\log n)/\sqrt{k})$ almost surely.
\end{cor}
We note that for Nadaraya-Watson estimate in Example 3, discussions in the next subsection suggest that the convergence behavior is very much the same under reasonable assumptions.

\textbf{Weakly dependent case} Here the convergence rate is determined by the interplay of $\psi$ and  $\{\gamma_m\}$ in a more complicated way. For example, qualitatively, the summability of $1/\psi(a/(2nv_{n1}\gamma_{m}))$ is easier to be satisfied the smaller is $\gamma_m$ (weaker dependence). Moreover, the choice of $x$ must take into account the trade off between the summability of $\exp\{-Ca^2/(aL+m^2c_{n2}^2+x)\}$ and the summability of $1/\psi(\sqrt{x/2}/(\gamma_1c_{n2}))$ (the former is an increasing function of $x$ while the latter is a decreasing function of $x$). Similarly, the choice of $m$ must take into account the trade off between summability of $(m/a)/\psi^{1-1/r}\left(L/(2Mmv_{n1})\right)$ and $1/\psi(a/(2nv_{n1}\gamma_{m}))$ (the former is an increasing function of $m$ while the latter is typically a decreasing function of $m$). Ignoring the complication of choosing $m$, when $\psi(x)=\psi_p(x)=\exp\{x^p\}-1$, the following corollary gives one possible situation where it is possible to set $a=mc_{n2}$ up to an extra logarithmic term.
\begin{cor}\label{cor:nw}
When $\psi=\psi_p, p\ge 1$, we have convergence rate $\|\hat{r}(x)-r(x)\|=O(b_n+H^\alpha+(\log n)^2mc_{n2}+(\gamma_1v_{n1})^{1/2})$ as long as $1/\psi(C(\log n)^2 m/(n\gamma_m))$ is summable for $C$ large enough.
\end{cor} 
\textit{Proof.} Take $x=C(\log n)^2c_{n2}^2$ ($C$ large enough) and $L_n=C(\log n)mc_{n2}$, the first expression in (*) is then satisfied if $a=C(\log n)^2mc_{n2}$. Moreover, $1/\psi(\sqrt{x/2}/(\gamma_1c_{n2}))\le 1/\psi(C\log n)$ is summable. Using the trivial inequalities $c_{n2}\ge v_{n1}$ and $v_{n1}\ge 1/n$, we get $m/a\le n$ and thus $(m/a)/\psi^{1-1/r}\left(L/(2Mmv_{n1})\right)\le n/\psi^{1-1/r}(C\log n)$ is summable. Finally, for the last sequence in (*), we have \[\sum_n1/\psi(a/(2nv_{n1}\gamma_{m_n}))\le \sum_n1/\psi(C(\log n)^2 m/(n\gamma_m))<\infty,\] by assumption in the statement of this corollary.

Finally, we note that in the above corollary, if $\gamma_m=e^{-Cm}$ for some $C>0$, then we can take $m\sim\log n$ so that all sequences in (*) are summable, and the rate of convergence is $(\log n)^3c_{n2}$.

\subsection{On the properties of $H$ and $k$ with random covariates}
In the previous subsection, we treat the predictor as fixed and the convergence rate depends on the sequence $\{X_i\}$. Here we study the behavior of some of the quantities that appeared in the rates when $X_i$ is a random stationary sequence in typical situations. Results obtained in this subsection can be combined with Theorem \ref{thm:main} to obtain more explicit convergence rates. The necessity of studying $H$ (for NN estimator) or $k$ (for Nadaraya-Watson estimator) is seen from Remark \ref{rmk:1} in the previous subsection.

When $X_i$ are random, we will make use of the important quantity $\varphi(h):=P(\{x':x'\in B(x,h)\})$ which is called the small ball probability. Its importance has been demonstrated in \cite{ferraty06} for functional kernel regression with scalar responses. In particular, the use of $\varphi(h)$ in a functional setting replaces the common assumption on the existence of a density for $X$ when $X$ belongs to some Euclidean space. It is easy to see that in the classical setting with mild assumption on the density of $X\in R^d$, we have $\varphi(h)\sim h^d$.

\textbf{Nearest neighbor estimate.} We only consider the simplest k-NN estimate as in Example 1 with $v_{ni}=1/k, i=1,\ldots,k$. Then in the convergence rates, $b_n=0$, $c_{n2}^2=\sum_iW_{ni}^2=1/k$ and $\max_iW_{ni}=1/k$. Thus only the quantity $H$ depends on $\{X_i\}$. If the sequence $\{X_i\}$ contains independent elements, one can show $H=O(\varphi(2k/n))$ almost surely as in the following proposition.

\begin{prop}\label{prop:1} Suppose $k/n\rightarrow 0$ and $k/\log n\rightarrow \infty$. Let $H$ be the distance from $x$ to its $k$-th nearest neighbor as defined in (\ref{eqn:H}), then $P(H> \varphi^{-1}(2k/n), i.o.)\rightarrow 0$, where $i.o.$ means ``infinitely often" and $\varphi^{-1}(x):=\inf\{h:\varphi(h)\ge x\}$. 
\end{prop}
\textit{Proof.} First we note that $\varphi$ is right-continuous and non-decreasing and thus $h=\phi^{-1}(x)$ implies $\phi(h)\ge x$.

 Denote $c=\varphi^{-1}(2k/n)$, $p=\varphi(c)$ and thus $np\ge 2k$. We have
\begin{eqnarray*}
&&P(H>\phi^{-1}(2k/n))\\
&=&P(\sum_iI\{X_i\in B(x,c)\}< k)\\
&=&P(\sum_iI\{X_i\in B(x,c)\}-np< k-np)\\
&\le&P(|\sum_iI\{X_i\in B(x,c)\}-np|\ge np/2)\\
&\le& 2\exp\{-\frac{1}{2}(np/2)^2/[np(1-p)+(np/6)]\}\\
&\le& 2\exp\{-Cnp\}\,,
\end{eqnarray*}
where we applied the Bernstein's inequality for Bernoulli random variables (see for example Appendix B in \cite{pollard84}). Then $P(H> \phi^{-1}(2k/n), i.o.)\rightarrow 0$ can be shown using Borel-Cantelli lemma noting that $k/\log n\rightarrow\infty$. $\Box$.

In \cite{ferraty06}, the authors distinguished two types of processes: the fractal type processes and the exponential type processes. The former is characterized by $\phi(h)\sim h^{\tau}$, for some $\tau>0$ and the latter characterized by $\phi(h)\sim \exp\{-(1/h^{\tau_1})\log(1/h^{\tau_2})\}, \tau_1>0,\tau_2\ge 0$. The fractal type processes are similar to finite dimensional covariates in many aspects, while for infinite dimensional case such as when the covariate curves belong to some smoothness class, exponential type processes are more typical. For example, the Brownian motion is of exponential type. The paper \cite{vaartzanten08} provides other more complicated Gaussian processes all of which are of exponential type. Combining Proposition \ref{prop:1} above with Corollary \ref{cor:knn}, we obtain the rates $O([\varphi^{-1}(2k/n)]^{\alpha}+(\log n)/\sqrt{k})$ for independent data. When the optimal $k$ is chosen, it is easy to see that for exponential type processes the convergence rates are logarithmic in the sample size, much slower than the classical finite-dimensional cases. Also note that this slow rate is largely determined by the term $[\varphi^{-1}(2k/n)]^{\alpha}$ which converges to zero logarithmically whether $k$ increases logarithmically or polynomially in $n$.

For weakly dependent sequence $\{X_i\}$, in particular assuming $\{X_i\}$ is $\psi-m-$approximable with $\|d(X_1,X_1^{(m)})\|=\beta_m, \sum_{m=1}^\infty\beta_m<\infty$ (a minor extension to Definition \ref{def:app} is needed here since $X_i\in\mathcal{F}$ which is not a normed space, thus we need to use $d(.,.)$ instead of $X_1-X_1^{(m)})$, we can show the following proposition whose proof is deferred to the next section. Note that although we used the same notation as before, $\psi$ here are different from that in Assumption 3 since here we are considering the predictor sequence instead of the noise sequence. 

\begin{prop}\label{prop:2}
Suppose for some $h>\varphi^{-1}(2k/n)$, there exists some sequence $1\le m\le n$ such that $k/n\rightarrow 0$, $k/(m\log n)\rightarrow\infty$ and $\sum_{n=1}^\infty n/\psi((h-\varphi^{-1}(2k/n))/\beta_m)<\infty$. Then we have $H\le h$ for $n$ large enough, almost surely.
\end{prop}
\textbf{Nadaraya-Watson estimate.} Here $W_{ni}=K(d(X_i,x)/H)/\sum_i K(d(X_i,x)/H)$ and we only consider the simple case where kernel function $K$ satisfies $  cI_{[-1,1]}\le K\le CI_{[-1,1]}$ for some $C>c>0$. Unlike k-NN estimate, here $H$ is predetermined. In Assumption 2, we let $k$ be the number of covariates inside the ball $B(x,H)$ and thus if $X_i$ is not one of the $k$ nearest neighbors of $x$, we have $W_{ni}=0$ and thus $b_n=\sum_{k+1}^nv_{ni}=0$ in the convergence rate in Theorem \ref{thm:main}. Since $H$ is predetermined in Nadaraya-Watson estimates, the only quantity in the convergence rates that depends on $X_i$ is $v_{n1}=\max_iW_{ni}$ and $c_{n2}=(\sum_iW_{ni}^2)^{1/2}$. Since $v_{n1}\le C/\sum_iK(d(X_i,x)/H)\le C/ck$, $v_{n1}\ge c/Ck$ as long as $k\ge 1$, and $c_{n2}\sim 1/\sqrt{k}$ which can be easily shown, we only need to study the asymptotic behavior of $k$, the number of predictors inside the ball $B(x,H)$.

With $\{X_i\}$ an independent sequence, we have
\begin{prop}\label{prop:3}
Suppose $H\rightarrow 0$, $n\varphi(H)/\log n\rightarrow\infty$, then $n\varphi(H)/2\le k\le 2n\varphi(H)$ for $n$ large enough, almost surely.
\end{prop}
On the other hand, for a $\psi-m-$approximable sequence $\{X_i\}$ with $\|d(X_1,X_1^{(m)})\|_\psi=\beta_m$, we have
\begin{prop}\label{prop:4}
Suppose $H''$ and $H'$ are two sequences with $H'<H<H''$ and there exists a sequence $1\le m\le n$ such that $n\varphi(H')/(m\log n)\rightarrow\infty$, $\sum_{n=1}^\infty n/\psi((H''-H)/\beta_m)<\infty$ and $\sum_{n=1}^\infty n/\psi((H-H')/\beta_m)<\infty$. Then we have $n\varphi(H')/2\le k\le 2n\varphi(H'')$ for $n$ large enough, almost surely.
\end{prop}
The proofs for these two propositions are very similar to those for Propositions \ref{prop:1} and \ref{prop:2}, and thus omitted.

\section{Proofs}
Based on two different representations of the nonparametric estimate in (\ref{eqn:rhat}) and (\ref{eqn:rhat2}), we decompose $\|\hat{r}(x)-r(x)\|$ into the bias term and the variance term,
\begin{equation}\label{bv}
\|\hat{r}(x)-r(x)\|\le\|\sum_i v_{ni}(r(X_{R_i})-r(x))\|+\|\sum_iW_{ni}\epsilon_i\|.
\end{equation}
The bias term is easier to deal with. In fact, 
\begin{eqnarray*}
\|\sum_i v_{ni}(r(X_{R_i})-r(x))\|&\le& 2B\sum_{i=k+1}^n v_{ni}+\|\sum_{i=1}^k v_{ni}(r(X_{R_i})-r(x))\|\\
&=&O(b_n+H^\alpha)\,,
\end{eqnarray*}
by Assumptions 1 and 2.

Now we deal with the variance term. Let $\eta_i=W_{ni}\epsilon_i$, $S_{n}=\sum_{i=1}^n\eta_i$ and the following arguments are conditional on $\{X_1,\ldots,X_n\}$ (in effect treating $W_{ni}$ as nonrandom weights). Following the idea of Section 6.3 in \cite{ledoux91}, we write $\|S_{n}\|-E\|S_n\|=\|\sum_{i=1}^n\eta_i\|-E\|\sum_{i=1}^n\eta_i\|=\sum_{i=1}^n e_i$, with $e_i=E[\|S_{n}\|\,|\mathcal{G}_i]-E[\|S_{n}\|\,|\mathcal{G}_{i-1}]$ where $\mathcal{G}_i$ is the $\sigma-$algebra generated by $\epsilon_1,\ldots,\epsilon_i$ ($\mathcal{G}_0$ is the trivial $\sigma-$algebra). It is easy to see that $\{e_i\}$ is a \textit{real-valued} martingale difference sequence which potentially enables us to use relevant exponential type inequalities. However, in general it seems at least not easy to obtain directly appropriate moment bounds for $d_i$ in order to apply, for example, Lemma 8.9 in \cite{geer00} (Bernstein's inequality for martingale differences), and thus we instead work with the quantity
\[
d_i=E[\|S_{n}\|\,|\mathcal{G}_i]-E[\|S_{n}\|\,|\mathcal{G}_{i-1}]-E[\|S_{n}-\eta_i-\cdots -\eta_{i+m-1}\|\,|\mathcal{G}_i]+E[\|S_{n}-\eta_i-\cdots -\eta_{i+m-1}\|\,|\mathcal{G}_{i-1}],
\]
where $m$ is same as that in the statement of the theorem and, as discussed in Remarks following the theorem, need to be chosen appropriately (as a side note, $m=1$ suffices for independent data in which case we actually have $d_i=e_i$). If $i+m-1>n$, the expression $S_n-\eta_i-\cdots-\eta_{i+m-1}$ is taken to mean $S_n-\eta_i-\cdots-\eta_n$. Obviously $d_i$ is still a martingale difference sequence. We denote $f_i=E[\|S_{n}-\eta_i-\cdots -\eta_{i+m-1}\|\,|\mathcal{G}_i]-E[\|S_{n}-\eta_i-\cdots -\eta_{i+m-1}\|\,|\mathcal{G}_{i-1}]$ and thus $e_i=d_i+f_i$.

Lemma \ref{lem:d} shows that 
\begin{equation}\label{eqn:d1}
|d_i|\le \sum_{j=i}^{i+m-1}W_{nj}E(\|\epsilon_j\||\mathcal{G}_i)+\sum_{j=i}^{i+m-1}W_{nj}E(\|\epsilon_j\||\mathcal{G}_{i-1}),
\end{equation}
and
\begin{equation}\label{eqn:d2}
E(d_i^2|\mathcal{G}_{i-1})\le m\sum_{j=i}^{i+m-1}W_{nj}^2E(\|\epsilon_j\|^2|\mathcal{G}_{i-1}).
\end{equation}
Lemma \ref{lem:f} shows that 
\begin{equation}\label{eqn:f}
P(\sum_if_i>a)\le 2/\psi\left(a/(2nv_{n1}\gamma_m)\right).
\end{equation}
Lemma \ref{lem:es} shows that 
\begin{equation}\label{eqn:es}
E\|S_n\|=O(c_{n2}+\sqrt{\gamma_1v_{n1}}).
\end{equation}

Aided by these results, we can bound the variance term $\|S_n\|$ in three steps.

\textit{Step 1: Let $d_i'=d_iI\{|d_i|\le L)$ for some $L>0$. We have $P(\sum_{i=1}^n (d_i'-E(d_i'|\mathcal{G}_{i-1}))>a)\le \exp\{-Ca^2/(aL+m^2c_{n2}^2+x)\}+1/\psi\left(\sqrt{x}/(\sqrt{2}\gamma_1c_{n2})\right), \forall a\ge 0, x\ge 0$.}

Let $\tilde{\psi}(x):=\psi(\sqrt{x})$. By Assumption 3, $\tilde{\psi}$ is convex and increasing and thus defines an Orlicz norm. Using (\ref{eqn:d2}), we have
\begin{eqnarray*}
&&\sum_{i=1}^n E(d_i^2|\mathcal{G}_{i-1})\\
&\le&m\sum_{i=1}^n \sum_{j=i}^{i+m-1}W_{nj}^2 E(\|\epsilon_j\|^2|\mathcal{G}_{i-1})\\
&=&m\sum_{i=1}^n \sum_{j=i}^{i+m-1}W_{nj}^2 E(\|\epsilon_j^{(j-i+1)}+\epsilon_j-\epsilon_j^{(j-i+1)}\|^2|\mathcal{G}_{i-1})\\
&\le& 2m^2E(\|\epsilon_1\|^2)\sum_{i=1}^nW_{ni}^2+2\sum_{i=1}^n \sum_{j=i}^{i+m-1}W_{nj}^2E(\|\epsilon_j-\epsilon_j^{(j-i+1)}\|^2|\mathcal{G}_{i-1}),
\end{eqnarray*}
where in the last line above we use that $\epsilon_j^{(j-i+1)}$ is independent of $\mathcal{G}_{i-1}$, and also use the inequality $\|\epsilon_j^{(j-i+1)}+\epsilon_j-\epsilon_j^{(j-i+1)}\|^2\le 2\|\epsilon_j^{(j-i+1)}\|^2+2\|\epsilon_j-\epsilon_j^{(j-i+1)}\|^2$ which follows from the parallelogram identity. Furthermore,
\begin{eqnarray*}
&&\|2\sum_{i=1}^n \sum_{j=i}^{i+m-1}W_{nj}^2E(\|\epsilon_j-\epsilon_j^{(j-i+1)}\|^2|\mathcal{G}_{i-1})\|_{\tilde{\psi}}\\
&\le&2\sum_{i=1}^n \sum_{j=i}^{i+m-1}W_{nj}^2 \left\|E(\|\epsilon_j-\epsilon_j^{(j-i+1)}\|^2|\mathcal{G}_{i-1})\right\|_{\tilde{\psi}}\\
&\le&2\sum_{i=1}^n \sum_{j=i}^{i+m-1}W_{nj}^2 \left\|\,\|\epsilon_j-\epsilon_j^{(j-i+1)}\|^2\,\right\|_{\tilde{\psi}}\\
&\le&2\sum_{i=1}^n \sum_{j=i}^{i+m-1}W_{nj}^2 \left(\|\epsilon_j-\epsilon_j^{(j-i+1)}\|_{{\psi}}\right)^{2}\\
&\le& 2\gamma_1^2 \sum_iW_{ni}^2,
\end{eqnarray*}
where we used Lemma \ref{lem:psi} (vii) for the second inequality above and Lemma \ref{lem:psi} (vi) for the third inequality above. Then, using Lemma \ref{lem:psi} (i), we have for any $x\ge 0$
\begin{equation}\label{eqn:R}
P(\sum_{i=1}^n E(d_i^2|\mathcal{G}_{i-1})>2m^2E(\|\epsilon_1\|^2)c_{n2}^2+x)\le 1/\psi(\sqrt{x}/(\sqrt{2}\gamma_1c_{n2})).
\end{equation}

Using $|d_i'-E(d_i'|\mathcal{G}_{i-1})|\le 2L$ and  $E[(d_i'-E(d_i'|\mathcal{G}_{i-1}))^2|\mathcal{G}_{i-1}]\le E(d_i'^2|\mathcal{G}_{i-1})\le E(d_i^2|\mathcal{G}_{i-1})$, we get $E(|d_i'-E(d_i'|\mathcal{G}_{i-1})|^k|\mathcal{G}_{i-1})\le (2L)^{k-2}E(d_i^2|\mathcal{G}_{i-1}), \forall k\ge 2$. Since $d_i'-E(d_i'|\mathcal{G}_{i-1}), i\le n$ is a martingale difference sequence, using Lemma 8.9 in \cite{geer00} (Bernstein's inequality for martingales) together with (\ref{eqn:R}), we obtain the desired bound as follows:
\begin{eqnarray*}
&&P(\sum_{i=1}^n(d_i'-E(d_i'|\mathcal{G}_{i-1}))>a)\\
&\le&P(\sum_{i=1}^n(d_i'-E(d_i'|\mathcal{G}_{i-1}))>a \mbox{ and } \sum_{i=1}^nE(d_i^2|\mathcal{G}_{i-1})\le 2m^2E(\|\epsilon_1\|^2)c_{n2}^2+x)\\
&&\;\;\;+P(\sum_{i=1}^nE(d_i^2|\mathcal{G}_{i-1})>2m^2E(\|\epsilon_1\|^2)c_{n2}^2+x)\\
&\le& \exp\{-Ca^2/(aL+m^2c_{n2}^2+x)\}+1/\psi\left(\sqrt{x}/(\sqrt{2}\gamma_1c_{n2})\right).
\end{eqnarray*}

\textit{Step 2: Let $d_i''=d_i-d_i'=d_iI\{|d_i|>L\}$. We have $P(\sum_i |d_i''-E(d_i''|\mathcal{G}_{i-1})|>a)\le Cm/\left(a\psi^{1-1/r}\left(\frac{L}{2Mmv_{n1}}\right)\right)$, where $M=\|\epsilon_1\|_\psi$.}

From (\ref{eqn:d1}), we have that
\begin{eqnarray*}
\|d_i\|_\psi&\le&2\sum_{j=i}^{i+m-1}W_{nj}\|\epsilon_j\|_\psi=2M\sum_{j=i}^{i+m-1}W_{nj},\\
\end{eqnarray*}
and thus using Lemma \ref{lem:psi} (i) and (v),
\[P(d_i>L)\le 1/\psi(\frac{L}{2M\sum_{j=i}^{i+m-1}W_{nj}}),\]
and
\[\|d_i\|_r\le C\|d_i\|_\psi\le C\sum_{j=i}^{i+m-1}W_{nj}.\]

Using H\"{o}lder's inequality, we have
\begin{eqnarray*}
&&E(|d_i''-E(d_i''|\mathcal{G}_{i-1})|)\\
&\le& 2E(|d_i''|)\\
&=&2E(|d_i|I\{|d_i|>L\})\\
&\le& 2\{E(|d_i|^r)\}^{1/r}P(|d_i|>L)^{1-1/r}\\
&\le& C(\sum_{j=i}^{i+m-1}W_{nj})/\psi^{1-1/r}\left(\frac{L}{2M\sum_{j=i}^{i+m-1}W_{nj}}\right)\\
&\le& C(\sum_{j=i}^{i+m-1}W_{nj})/\psi^{1-1/r}\left(\frac{L}{2Mmv_{n1}}\right),
\end{eqnarray*}
and thus, using Markov's inequality, we have $P(\sum_i |d_i''-E(d_i''|\mathcal{G}_{i-1})|>a)\le E[\sum_i |d_i''-E(d_i''|\mathcal{G}_{i-1})|]/a\le Cm/\left(a\psi^{1-1/r}\left(\frac{L}{2Mmv_{n1}}\right)\right)$.

\textit{Step 3: Finally, we demonstrate the bound for the variance term in (\ref{bv}).}

Using $E(d_i|\mathcal{G}_{i-1})=E(d_i'|\mathcal{G}_{i-1})+E(d_i''|\mathcal{G}_{i-1})=0$, we have that $d_i=d_i'-E(d_i'|\mathcal{G}_{i-1})+(d_i''-E(d_i''|\mathcal{G}_{i-1}))$ and then
\begin{eqnarray*}
&&P(\|S_{n}\|-E\|S_{n}\|>3a)\\
&=&P(\sum_i d_i+f_i>3a)\\
&\le&P(\sum_i d_i>2a)+P(\sum_if_i>a)\\
&\le&P(\sum_i (d_i'-E(d_i'|\mathcal{G}_{i-1}))>a)+P(\sum_i (d_i''-E(d_i''|\mathcal{G}_{i-1}))>a)+P(\sum_if_i>a)\\
&\le&\exp\{-Ca^2/(aL+m^2c_{n2}^2+x)\}+1/\psi\left(\sqrt{x}/(\sqrt{2}\gamma_1c_{n2})\right)\\
&&+C(m/a)/\psi^{1-1/r}\left(\frac{L}{2Mmv_{n1}}\right)+2/\psi\left(a/(2nv_{n1}\gamma_m)\right),
\end{eqnarray*}
by the previous two steps and (\ref{eqn:f}). The above expression is summable by assumption of the theorem, and an application of the Borel-Cantelli Lemma leads to $\|S_{n}\|-E\|S_{n}\|=O(a)$. Combining this with (\ref{eqn:es}), the variance term is thus $\|S_{n}\|=O(a+c_{n2}+(\gamma_1v_{n1})^{1/2})$. As noted in Remark 1 following the theorem, the term $c_{n2}$ can be omitted since we always have $c_{n2}=O(a)$.$\Box$

\begin{lem}\label{lem:d} Using the notation in the proof of Theorem \ref{thm:main}, we have 
\begin{equation*}
|d_i|\le \sum_{j=i}^{i+m-1}E(\|\eta_j\||\mathcal{G}_i)+\sum_{j=i}^{i+m-1}E(\|\eta_j\||\mathcal{G}_{i-1}),
\end{equation*}
\begin{equation*}
E(d_i^2|\mathcal{G}_{i-1})\le m\sum_{j=i}^{i+m-1}E(\|\eta_j\|^2|\mathcal{G}_{i-1}).
\end{equation*}
\end{lem}
\textit{Proof.}
Since $d_i=E(\|S_n\|-\|S_n-\sum_{j=i}^{i+m-1}\eta_j\||\mathcal{G}_i)-E(\|S_n\|-\|S_n-\sum_{j=i}^{i+m-1}\eta_j\||\mathcal{G}_{i-1})$, the first equation is obvious.

Denote $\xi_i=E(\|S_n\|-\|S_n-\sum_{j=i}^{i+m-1}\eta_j\||\mathcal{G}_i)$, then $d_i=\xi_i-E(\xi_i|\mathcal{G}_{i-1})$. Using the interpretation of $E(\xi_i|\mathcal{G}_{i-1})$ as the projection of $\xi_i$, we have
\[E(d_i^2|\mathcal{G}_{i-1})\le E(\xi_i^2|\mathcal{G}_{i-1})\le m\sum_{j=i}^{i+m-1}E(\|\eta_j\|^2|\mathcal{G}_{i-1}),\]
proving the second equation.
$\Box$

\begin{lem}\label{lem:f} For $f_i=E[\|S_{n}-\eta_i-\cdots -\eta_{i+m-1}\|\,|\mathcal{G}_i]-E[\|S_{n}-\eta_i-\cdots -\eta_{i+m-1}\|\,|\mathcal{G}_{i-1}]$ as in the proof of Theorem \ref{thm:main}, we have
\[P(\sum_if_i>a)\le 2/\psi\left(a/(2n(\max_{1\le j\le n}W_{nj})\gamma_m)\right).\]
\end{lem}
\textit{Proof.} By the definition of $\psi-m-$approximability for sequence $\epsilon_i$, we have that
\begin{eqnarray*}
&&E[\|W_{n1}\epsilon_1+\cdots+W_{n i-1}\epsilon_{i-1}+W_{n i+m}\epsilon_{i+m}^{(m)}+\cdots+W_{nn}\epsilon_n^{(n-i)}\|\,|\mathcal{G}_i]\\
&=&E[\|W_{n1}\epsilon_1+\cdots+W_{n i-1}\epsilon_{i-1}+W_{n i+m}\epsilon_{i+m}^{(m)}+\cdots+W_{nn}\epsilon_n^{(n-i)}\|\,|\mathcal{G}_{i-1}].
\end{eqnarray*}
Thus $f_i=f_i'-f_i''$ where 
\begin{eqnarray*}
f_i'&=&E[\|W_{n1}\epsilon_1+\cdots+W_{n i-1}\epsilon_{i-1}+W_{n i+m}\epsilon_{i+m}+\cdots+W_{nn}\epsilon_n\|\,|\mathcal{G}_i]\\
&&-E[\|W_{n1}\epsilon_1+\cdots+W_{n i-1}\epsilon_{i-1}+W_{n i+m}\epsilon_{i+m}^{(m)}+\cdots+W_{nn}\epsilon_n^{(n-i)}\|\,|\mathcal{G}_i],\\
f_i''&=&E[\|W_{n1}\epsilon_1+\cdots+W_{n i-1}\epsilon_{i-1}+W_{n i+m}\epsilon_{i+m}+\cdots+W_{nn}\epsilon_n\|\,|\mathcal{G}_{i-1}]\\
&&-E[\|W_{n1}\epsilon_1+\cdots+W_{n i-1}\epsilon_{i-1}+W_{n i+m}\epsilon_{i+m}^{(m)}+\cdots+W_{nn}\epsilon_n^{(n-i)}\|\,|\mathcal{G}_{i-1}].
\end{eqnarray*}
Since $|f_i'|\le E(\|W_{ni+m}\epsilon_{i+m}-W_{ni+m}\epsilon_{i+m}^{(m)}\|+\cdots+\|W_{nn}\epsilon_{n}-W_{nn}\epsilon_{n}^{(n-i)}\|\,|\mathcal{G}_i)$, using Lemma \ref{lem:psi} (vii), we have $\|f_i'\|_\psi\le (\max_{1\le j\le n}W_{nj})\gamma_m$ and thus $\|\sum_{i=1}^n f_i'\|_\psi\le n(\max_{1\le j\le n}W_{nj})\gamma_m$. Using Lemma \ref{lem:psi} (i) we get
\[P(\sum_if_i'>a/2)\le 1/\psi\left(a/(2n(\max_{1\le j\le n}W_{nj})\gamma_m)\right).\]
By exactly the same arguments
\[P(\sum_if_i''>a/2)\le 1/\psi\left(a/(2n(\max_{1\le j\le n}W_{nj})\gamma_m)\right),\]
and the Lemma is proved by combining the above two displayed equations.
$\Box$

\begin{lem}\label{lem:es} Let $S_n=\sum_{i=1}^n\eta_i=\sum_{i=1}^n W_{ni}\epsilon_i$ as in the proof of Theorem \ref{thm:main}, we have $E\|S_n\|=O\left(c_{n2}+\sqrt{\gamma_1\max_iW_{ni}}\right)$.
\end{lem}
\textit{Proof.} We have 
\begin{eqnarray*}
&&(E\|S_n\|)^2\\
&=&(E\|\sum_iW_{ni}\epsilon_i\|)^2\\
&\le& E(\|\sum_iW_{ni}\epsilon_i\|^2)\\
&=&E\sum_i\sum_jW_{ni}W_{nj}\langle\epsilon_i,\epsilon_j\rangle\\
&=&\sum_i W_{ni}^2E\|\epsilon_i\|^2+2\sum_{i=1}^n\sum_{j=i+1}^nW_{ni}W_{nj}E\langle\epsilon_i,\epsilon_j\rangle\\
&=&O(c_{n2}^2)+2\sum_{i=1}^n\sum_{j=i+1}^nW_{ni}W_{nj}E\langle\epsilon_i,\epsilon_j-\epsilon_j^{(j-i)}\rangle\\
&\le&O(c_{n2}^2)+2\sum_{i=1}^n\sum_{j=i+1}^nW_{ni}W_{nj}E(\|\epsilon_i\|\,\|\epsilon_j-\epsilon_j^{(j-i)}\|)\\
&\le&O(c_{n2}^2)+2\sum_{i=1}^n\sum_{j=i+1}^nW_{ni}W_{nj}(E(\|\epsilon_i\|^2)^{1/2} (E\|\epsilon_j-\epsilon_j^{(j-i)}\|^2)^{1/2}\\
&\le&O(c_{n2}^2)+C\sum_{i=1}^n\sum_{j=i+1}^nW_{ni}W_{nj}\|\epsilon_i\|_\psi \|\epsilon_j-\epsilon_j^{(j-i)}\|_\psi,
\end{eqnarray*}
where we used that $\epsilon_j$ and $\epsilon_j^{(j-i)}$ are independent, and Assumption 3 on $\psi$. Finally, we see that $\sum_{i=1}^n\sum_{j=i+1}^nW_{ni}W_{nj}\|\epsilon_i\|_\psi \|\epsilon_j-\epsilon_j^{(j-i)}\|_\psi \le (\max_jW_{nj})(\sum_iW_{ni})\sum_{m=1}^\infty \|\epsilon_1-\epsilon_1^{(j-i)}\|_\psi=O\left(\gamma_1\max_iW_{ni}\right).$
$\Box$

\textit{Proof of Proposition \ref{prop:2}.}
Consider the approximating sequence $(X_1^{(m)},X_2^{(m)},\ldots)$. Define the zero mean random variables $Y_i^{(m)}=I\{X_i^{(m)}\in B(x,c)\}-\tilde{p}$ where $c=\varphi^{-1}(2k/n)$ and $\tilde{p}=\varphi(c)=2k/n$. Divide the sequence $(Y_1^{(m)},Y_2^{(m)},\ldots)$ into $m$ groups (we assume $n/m$ is an integer for simplicity in presentation without loss of generality) as follows:
\[
\begin{array}{rl}
\mbox{group 1:} & Y_1^{(m)}, Y_{1+m}^{(m)}, Y_{1+2m}^{(m)},\ldots, Y_{1+(n/m-1)m}^{(m)},\\
\mbox{group 2:} & Y_2^{(m)}, Y_{2+m}^{(m)}, Y_{2+2m}^{(m)},\ldots, Y_{2+(n/m-1)m}^{(m)},\\
\vdots & \vdots\\
\mbox{group $m$:} & Y_m^{(m)}, Y_{2m}^{(m)}, Y_{3m}^{(m)},\ldots, Y_{n}^{(m)}.\\
\end{array}
\]
Because of the construction, the random variables within one group are independent of each other. Let $Z_i, i=1,\ldots m$ be the sum of random variables within each group. Using Bernstein's inequality, we have 
\[P(|Z_i|>x)\le 2\exp\{-\frac{1}{2}x^2/(n\tilde{p}/m+x/3)\},\]
and thus 
\begin{eqnarray}
&&P(\sum_{i=1}^nI\{X_i^{(m)}\in B(x,c)\}\le k)\nonumber\\
&\le& P(|\sum_{i=1}^n(I\{X^{(m)}_i\in B(x,c)\}-\tilde{p})|\ge n\tilde{p}-k)\nonumber\\
&=&P(|\sum_{i=1}^m Z_i|\ge n\tilde{p}-k)\nonumber\\
&\le & mP(|Z_1|>(n\tilde{p}-k)/m)\nonumber\\
&\le &2m \exp\{-\frac{1}{2}(\frac{n\tilde{p}-k}{m})^2/(n\tilde{p}/m+(n\tilde{p}-k)/(3m))\nonumber\\
&=& 2m\, \exp\{-(3/14)k/m\}.\label{eqn:H1}
\end{eqnarray}


We also have that 
\begin{eqnarray}
&&P(H>h)\le P(\sum_{i=1}^nI\{X_i\in B(x,h)\}\le k)\nonumber\\
&\le& P(\sum_{i=1}^nI\{X_i^{(m)}\in B(x,c)\}\le k)\nonumber\\
&&\;\;\;\;+P(\exists i, s.t.\; X_i^{(m)}\in B(x,c) \mbox{ and }  X_i\not\in B(x,h))\nonumber\\
&\le&P(\sum_{i=1}^nI\{X_i^{(m)}\in B(x,c)\}\le k)+P(\exists i, s.t. \; d(X_i^{(m)},X_i)>h-c)\nonumber\\
&\le&P(\sum_{i=1}^nI\{X_i^{(m)}\in B(x,c)\}\le k)+n/\psi((h-c)/\beta_m)\label{eqn:ac}\\
&\le& 2m\, \exp\{-(3/14)k/m\}+n/\psi((h-c)/\beta_m),\label{eqn:H2}
\end{eqnarray}
where we used Lemma \ref{lem:psi} (i) in (\ref{eqn:ac}) and used (\ref{eqn:H1}) in (\ref{eqn:H2}). The lemma follows from the Borel-Cantelli Lemma.
$\Box$.

\bibliographystyle{plain}
\bibliography{papers.txt,books.txt}

\end{document}